\documentclass[sigconf]{acmart}

\usepackage{amsmath,amsfonts}
\usepackage{graphicx}
\usepackage{algorithmic}
\usepackage{subcaption}
\usepackage{comment}
\usepackage{graphicx}
\usepackage{textcomp}
\usepackage[linesnumbered,ruled,vlined]{algorithm2e}
\usepackage{listings}
\usepackage{xspace}
\usepackage{xcolor}
\usepackage{cleveref}
\usepackage{multirow}

\def\BibTeX{{\rm B\kern-.05em{\sc i\kern-.025em b}\kern-.08em
    T\kern-.1667em\lower.7ex\hbox{E}\kern-.125emX}}
\definecolor{Yellow2}{HTML}{C6BD2F}

\newcommand\norm[1]{\left\lVert#1\right\rVert}

\newcommand{\redtext}[1]{{\color[HTML]{d73027}{#1}}}

\newcommand{\darkbluetext}[1]{{\color[HTML]{4575b4}{#1}}}
\newcommand{\orangetext}[1]{{\color[HTML]{fc8d59}{#1}}}

\newcommand{\parsec}{PaRSEC\xspace}
\newcommand{\hatrix}{$\mathcal{H}$ATRIX\xspace}
\newcommand{\lorapo}{LORAPO\xspace}
\newcommand{\strumpack}{STRUMPACK\xspace}

\AtBeginDocument{%
  \providecommand\BibTeX{{%
    \normalfont B\kern-0.5em{\scshape i\kern-0.25em b}\kern-0.8em\TeX}}}

\copyrightyear{2023} 
\acmYear{2023} 
\setcopyright{acmlicensed}\acmConference[ICPP 2023]{52nd International Conference on Parallel Processing}{August 7--10, 2023}{Salt Lake City, UT, USA}
\acmBooktitle{52nd International Conference on Parallel Processing (ICPP 2023), August 7--10, 2023, Salt Lake City, UT, USA}
\acmPrice{15.00}
\acmDOI{10.1145/3605573.3605606}
\acmISBN{979-8-4007-0843-5/23/08}

\begin{document}

\title{$O(N)$ distributed direct factorization of structured dense matrices using runtime systems.}

\author{Sameer Deshmukh}
\email{sameer.deshmukh@rio.gsic.titech.ac.jp}
\affiliation{%
  \institution{School of Computing, Tokyo Institute of Technology}
  \city{Tokyo}
  \country{Japan}
  \postcode{152-8550}
}
\orcid{0000-0002-5615-1399}

\author{Qinxiang Ma}
\email{ma@rio.gsic.titech.ac.jp}
\affiliation{%
  \institution{School of Computing, Tokyo Institute of Technology}
  \city{Tokyo}
  \country{Japan}
}
\orcid{0000-0002-7329-1756}

\author{Rio Yokota}
\email{rioyokota@gsic.titech.ac.jp}
\affiliation{%
  \institution{Global Scientific Information and Computing Center, Tokyo Institute of Technology}
  \city{Tokyo}
  \country{Japan}
  \postcode{152-8550}
}
\orcid{0000-0001-7573-7873}

\author{George Bosilca}
\email{bosilca@icl.utk.edu}
\affiliation{%
  \institution{Innovative Computing Laboratory, University of Tennessee at Knoxville}
  \city{Knoxville}
  \country{USA}
}
\orcid{0000-0003-2411-8495}


\begin{abstract}
Structured dense matrices result from boundary integral problems in electrostatics and geostatistics, and also Schur complements in sparse preconditioners such as multi-frontal methods. Exploiting the structure of such matrices can reduce the time for dense direct factorization from $O(N^3)$ to $O(N)$. The Hierarchically Semi-Separable (HSS) matrix is one such low rank matrix format that can be factorized using a Cholesky-like algorithm called ULV factorization. The HSS-ULV algorithm is highly parallel because it removes the dependency on trailing sub-matrices at each HSS level. However, a key merge step that links two successive HSS levels remains a challenge for efficient parallelization. In this paper, we use an asynchronous runtime system \parsec with the HSS-ULV algorithm. We compare our work with \strumpack and \lorapo, both state-of-the-art implementations of dense direct low rank factorization, and achieve up to 2x better factorization time for matrices arising from a diverse set of applications on up to 128 nodes of Fugaku for similar or better accuracy for all the problems that we survey.

\end{abstract}




\maketitle

\section{Introduction}
\label{sec:introduction}

Various scientific problems in fluid dynamics, structural mechanics, and electromagnetics are governed by partial differential equations (PDEs), which result in a sparse matrix when discretized with finite difference and finite element methods. The boundary element method (BEM) has an advantage over such volume discretization methods, since it only discretizes the boundary and the number of elements can be drastically reduced.
However, BEM requires the solution of a dense linear system, which has cubic complexity if solved directly.
For very large dense matrices, the Cholesky factorization can be computed using distributed memory implementations from SCALAPACK, DPLASMA~\cite{bosilca2011}, SLATE~\cite{gatesSLATEDesignModern2019} or Elemental~\cite{poulson2013}.
Although modern parallel computer architectures offer significant speedups due to the use of multiple threads, the cubic complexity of dense factorization remains prohibitive for large matrices.

\begin{table*}[]
\begin{tabular}{|c|c|c|c|c|c|}
  \hline
  \textbf{Library} & \textbf{Format} & \textbf{Algorithm} & \textbf{Compute complexity} & \textbf{Distributed paradigm} & \textbf{Comm. complexity} \\
  \hline
  DPLASMA~\cite{bosilca2011}                              & Dense & Tile Cholesky & $O(N^3)$ & Asynchronous  & $O(N^3)$ \\
  \hline
  SLATE~\cite{gatesSLATEDesignModern2019}                 & Dense & Panel Cholesky & $O(N^3)$ & Fork-join  & $O(N^3)$ \\
  \hline
  LORAPO~\cite{cao2022a}                                   & BLR & Tile Cholesky & $O(N^2)$ & Asynchronous  & $O(N^3)$ \\
  \hline
  $\mathcal{H}$-LU~\cite{augonnetHierarchicalFastDirect2019}          & $\mathcal{H}$-matrix & $\mathcal{H}$-LU & $O(Nlog(N))$ & Asynchronous & $O(Nlog(N))$ \\
  \hline
  STRUMPACK~\cite{rouetDistributedmemoryPackageDense2015} & HSS & ULV & $O(N)$ & Fork-join & $O(N^2)$ \\
  \hline
  Ma et. al.~\cite{ma2022a}                                     & $\mathcal{H}^2$-matrix & Modified ULV & $O(N)$ & Fork-join & $O(N)$ \\
  \hline
  \hatrix-DTD & HSS & ULV & $O(N)$ & Asynchronous & $O(N)$\\
  \hline
\end{tabular}
\caption{Comparison of dense direct factorization methods depending on the matrix format, factorization algorithm and distributed programming paradigm.}
\label{tab:formats-and-algorithms}
\end{table*}

Table~\ref{tab:formats-and-algorithms} shows some of the state-of-the-art implementations of distributed dense direct factorization using various matrix formats and algorithms. Going from top to bottom, we can see libraries such as DPLASMA~\cite{bosilca2011} and SLATE~\cite{gatesSLATEDesignModern2019} (and SLATE's predecessor SCALAPACK) making use of standard dense matrix formats and their associated algorithms that do not make use of low rank approximation. Although SLATE and DPLASMA use the same dense Cholesky factorization algorithm, DPLASMA makes use of asynchronous distributed execution whereas SLATE uses fork-join parallelism. The remaining libraries all make use of low rank representations of the dense matrix. This is achieved by compressing the off-diagonal blocks in the dense matrix that correspond to far interactions in the physical domain. The compression can be performed using a suitable algorithm such as Randomized Singular Value Decomposition (RSVD) or Adaptive Cross Approximation (ACA)~\cite{rjasanowAdaptiveCrossApproximation2002}. Hierarchical matrices exploit this property to reduce the cost of computation from $O(N^3)$ to almost $O(N^2)$ or even $O(N)$. The complexity of the factorization is determined by the format of the hierarchical matrix and the algorithm used for the factorization. Various formats such as BLR~\cite{amestoyImprovingMultifrontalMethods2015}, BLR$^2$~\cite{ashcraftBlockLowRankMatrices2021}, HODLR~\cite{ambikasaranNlogNFastDirect2013}, $\mathcal{H}$-matrix, HSS~\cite{chandrasekaran_ulv_2006} and $\mathcal{H}^2$-matrix~\cite{bormMatrixArithmeticsLinear2006} have been proposed, varying by conditions of admissibility and the use of nested basis.

The BLR format used by \lorapo~\cite{cao2022a} subdivides the dense matrix into blocks of uniform size and approximates each block individually. The BLR format reduces the time complexity of the tile Cholesky factorization to $O(N^2)$. The use of an asynchronous run time system such as \parsec allows \lorapo~\cite{cao2022a} to prioritize the execution of the critical path of the tile Cholesky factorization and resolve off-diagonal dependencies asynchronously. Large, adjacent low rank blocks of the BLR format can be combined to form the multi-level $\mathcal{H}$-matrix format. The tile LU (or Cholesky) factorization can then be extended to the $\mathcal{H}$-LU (or $\mathcal{H}$-Cholesky) algorithm which costs $O(Nlog(N))$. The use of an asynchronous runtime system with $\mathcal{H}$-LU has been shown to achieve good strong scaling for distributed computation since this allows for greater parallelism between the recursive blocks of the $\mathcal{H}$-matrix..

The use of nested basis in formats such as BLR$^2$, HSS and $\mathcal{H}^2$-matrix can be combined with the ULV factorization~\cite{chandrasekaran_ulv_2006} to further reduce the time complexity of factorization to close to $O(N)$. The ULV factorization of the HSS matrix (HSS-ULV) exploits the nested basis property to remove the dependency on the off-diagonal blocks during the Cholesky factorization. This means that there is no longer a need to perform the triangular solve and trailing sub-matrix update, which leads to an embarrassingly parallel factorization of successive levels of the HSS matrix. Therefore, dependencies only exist between the levels. Distributed HSS-ULV factorization has been implemented by \strumpack~\cite{rouetDistributedmemoryPackageDense2015} using the fork-join programming paradigm as a result of relying on SCALAPACK for the computation. \strumpack~\cite{rouetDistributedmemoryPackageDense2015} distributes each block of the HSS matrix with a block cyclic distribution and relies on collective communication to shuffle data. Ma et. al.~\cite{ma2022a} extend the ULV factorization to the $\mathcal{H}^2$-matrix, by modifying the ULV factorization to precompute the fill-ins before the factorization for the $\mathcal{H}^2$-matrix. Even though the $\mathcal{H}^2$-matrix has off-diagonal dense blocks, the method from Ma et. al.~\cite{ma2022a} is able to achieve embarrassingly parallel factorization of each level by performing the factorization twice - once for precomputing the fill-ins and then for the actual factorization. Although this method is highly parallel, the fact that it factorizes twice results in a large overhead.

In this paper, we propose \hatrix-DTD -- an implementation of the HSS-ULV factorization with the \parsec runtime system. \hatrix-DTD makes use of the HSS-ULV algorithm that can factorize matrices arising from a variety of Green's functions with comparable accuracy to \lorapo and \strumpack. 
We choose these codes as a reference because they share certain traits with our code, besides the fact that they are the most popular libraries in this field.
Similar to our code, \lorapo uses the \parsec runtime system, but the matrix structure is BLR.
Conversely, \strumpack uses the HSS structure like \hatrix-DTD, but uses a bulk-synchronous model for parallelism instead of a runtime system.
By comparing with these two references, we can isolate the effect of choosing HSS over BLR from the effect of using a runtime system over a bulk synchronous approach.
The HSS-ULV factorization is computationally cheaper than the BLR-Cholesky factorization, and hence \hatrix-DTD can outperform \lorapo while making use of the same runtime system. The use of an asynchronous runtime system such as \parsec for handling the communication and dependencies between successive levels in the HSS-ULV leads to better overlap of communication and computation, and hence \hatrix-DTD can outperform \strumpack that makes use of a similar HSS-ULV algorithm. As a result, we experimentally prove two key assumptions about the factorization algorithms and distributed memory implementation of low rank matrix formats in this paper:
\begin{enumerate}
\item The use of the HSS matrix format and HSS-ULV algorithm has lower computational complexity than the BLR-tile Cholesky algorithm implemented by \lorapo~\cite{cao2022a}.
\item The use of an asynchronous runtime system such as \parsec leads to a lower overhead of communication than the fork-join parallelism implemented by \strumpack~\cite{rouetDistributedmemoryPackageDense2015}.
\end{enumerate}

We experimentally show that \hatrix-DTD can outperform both \lorapo and \strumpack~\cite{rouetDistributedmemoryPackageDense2015} over a large number of nodes. \hatrix-DTD shows weak scaling efficiency and achieves up to 2x faster time of factorization on up to 128 nodes for a variety of Green's functions. We first introduce our notation and construction of HSS matrices, and then elaborate on the HSS-ULV algorithm. We then describe our implementation of the HSS-ULV using the \parsec runtime system. We then demonstrate with rigorous experimental evidence the performance of \hatrix-DTD against \strumpack~\cite{rouetDistributedmemoryPackageDense2015} and \lorapo. 

\section{Construction and notation}
\label{sec:construction-and-notation}

\begin{figure}
    \centering
    \includegraphics[width=0.9\linewidth]{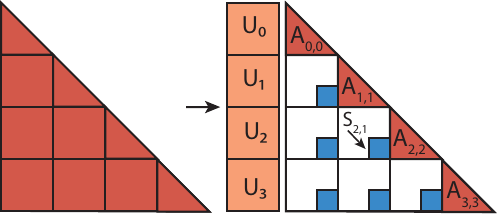}
    \caption{Construction of an SPD BLR$^2$ matrix from a block dense matrix.}
    \label{fig:blr2-construction}
\end{figure}

A symmetric positive definite BLR$^2$ matrix can be constructed from a block dense matrix as shown in Fig.~\ref{fig:blr2-construction}. A single block of this matrix at the index $(row,column)$ is denoted by \redtext{$A_{row,column}$}. We use a single shared bases denoted by \orangetext{$U_{row}$} to denote the bases of the admissible blocks on $row$. For example, the block \redtext{$A_{2,1}$} in Fig.~\ref{fig:blr2-construction} is denoted as
\begin{align}
    \redtext{A_{2,1}} \gets \orangetext{U_{2}} \cdot \darkbluetext{S_{2,1}} \cdot \orangetext{U_{1}^T}
\end{align}
where \darkbluetext{$S_{2,1}$} denotes the skeleton block shown in green.

The shared basis for each column is generated by concatenating the admissible blocks in the column, denoted by $A_{+,0}$. The shared basis can then be computed by computing a pivoted QR factorization
\begin{align}
    \orangetext{
        \begin{bmatrix}
            U^{S}_0 & U^R_{0}
        \end{bmatrix}
    } \gets QR(A_{+,0}^T)
\end{align}
where the $S$ and $R$ superscripts denote the skeleton part and redundant part of the basis, respectively. In order to make it convenient to represent the \orangetext{U}\darkbluetext{L}\orangetext{V} factorization, we permute the skeleton and redundant parts as shown in Eq.~(\ref{eq:complement-bases-ui}).

\begin{equation}
    \orangetext{U_{i}} = \orangetext{
        \begin{bmatrix}
            U^{R}_{i} & U^{S}_{i}
        \end{bmatrix}
    }
    \label{eq:complement-bases-ui}
\end{equation}

Any dense block of the BLR$^2$ matrix in Fig.~\ref{fig:blr2-construction} $A_{i,j}$ can be represented as Eq.~(\ref{eq:dense-block}).
\begin{align}
    \redtext{A_{i,j}} = \orangetext{\begin{bmatrix}
        U^R_i & U^S_i
    \end{bmatrix}} \cdot
    \darkbluetext{
        \begin{bmatrix}
            S^{RR}_{i,j} & S^{SR}_{i,j}\\
            S^{RS}_{i,j} & S^{SS}_{i,j}
        \end{bmatrix}
    }\cdot
    \orangetext{
        \begin{bmatrix}
            {U^R_j}^T \\
            {U^S_j}^T
        \end{bmatrix}
    }
    \label{eq:dense-block}
\end{align}
Similarly, any low rank block $A_{i,j}$ can be represented as Eq.~(\ref{eq:low-block}).
\begin{align}
    \redtext{A_{i,j}} =
        \orangetext{
            \begin{bmatrix}
                U^R_i & U^S_i
            \end{bmatrix}
    } \cdot
    \darkbluetext{\begin{bmatrix}
        0 & 0\\
        0 & S^{SS}_{i,j}
    \end{bmatrix}} \cdot
    \orangetext{\begin{bmatrix}
        {U^R_j}^T \\
        {U^S_j}^T
    \end{bmatrix}}
    \label{eq:low-block}
\end{align}

Low rank matrix formats that have dense blocks in their off diagonals are termed as being strongly admissible, and those that have dense blocks only on the diagonal are termed as weakly admissible. The notion of the weakly admissibility BLR$^2$ matrix described above can be extended to the HSS matrix. The HSS matrix introduces multiple levels in the matrix by sharing the basis between levels. This should not be confused with the recursive hierarchical structure of the HODLR~\cite{ambikasaranNlogNFastDirect2013} matrix, which does not share the basis but instead uses recursive low rank blocks in the off-diagonals.

\begin{figure}
    \centering
    \includegraphics[width=0.7\linewidth]{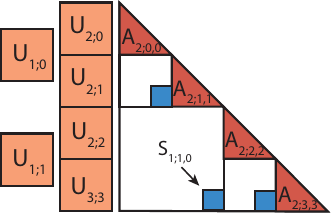}
    \caption{HSS construction and notation.}
    \label{fig:hss-construction}
\end{figure}

We introduce the notion of $level$ in order to represent the blocks of the HSS matrix at various levels. At the leaf level, the dense block \redtext{$A_{i,j}$} of the BLR$^2$ matrix can be represented as \redtext{$A_{level;i,j}$} for the HSS matrix. The shared basis at the leaf level also use a similar notation and are denoted by \orangetext{$U_{level;i}$}. Fig.~\ref{fig:hss-construction} introduces the notation and construction of a 2-level HSS matrix from the BLR$^2$ matrix shown in Fig.~\ref{fig:blr2-construction}. As an example, the block \redtext{$A_{1;1,0}$} can be represented with the nested basis as shown in Eq.~(\ref{eq:hss-block-notation}).

\begin{equation}
    \redtext{A_{1;1,0}} =
    \orangetext{\begin{bmatrix}
        U_{2;2} & 0 \\
        0 & U_{2;3}
    \end{bmatrix}} \cdot
    \orangetext{U_{1;1}} \cdot
    \darkbluetext{\begin{bmatrix}
        0 & 0 \\
        0 & S^{SS}_{1;1,0}
    \end{bmatrix}} \cdot
    \orangetext{U_{1;0}^T} \cdot
    \orangetext{\begin{bmatrix}
        U_{2;0}^T & 0 \\
        0 & U_{2;1}^T
    \end{bmatrix}}
    \label{eq:hss-block-notation}
\end{equation}

\section{The ULV factorization algorithm}
\label{sec:ulv-factorization-algorithm}
BLR$^2$ and HSS matrices can be factorized with the ULV factorization. The ULV factorization can be thought of as a modified Cholesky factorization where the $L$ represents the lower triangular dense blocks and the $U$ and $V$ represent the bases. The ULV works on the principle of nullifying the low rank off-diagonal blocks by multiplying the row and column with the shared bases. This means that there is no need to perform the triangular solve and trailing sub-matrix updates for the admissible blocks.

Multiplication of the dense block shown in Eq.~(\ref{eq:dense-block}) with its respective row and column basis from Eq.~(\ref{eq:complement-bases-ui}) leads to
\begin{align}
    \darkbluetext{\begin{bmatrix}
        S^{RR}_{i,j} & S^{SR}_{i,j}\\
        S^{RS}_{i,j} & S^{SS}_{i,j}
    \end{bmatrix}} =
    \orangetext{\begin{bmatrix}
        {U^{R}}^T_i \\ {U^{S}}^T_i
    \end{bmatrix}} \cdot
    \redtext{A_{i,j}}
    \orangetext{\begin{bmatrix}
        {U^{R}}_j & {U^{S}}_j
    \end{bmatrix}}
    \label{eq:dense-split}
\end{align}
Likewise, multiplication of the low rank block from Eq.~(\ref{eq:low-block}) leads to Eq.~(\ref{eq:low-rank-split}).
\begin{align}
    \darkbluetext{\begin{bmatrix}
        0 & 0\\
        0 & S^{SS}_{i,j}
    \end{bmatrix}} =
    \orangetext{\begin{bmatrix}
        {U^{R}}^T_i \\ {U^{S}}^T_i
    \end{bmatrix}} \cdot
    \redtext{A_{i,j}}
    \orangetext{\begin{bmatrix}
        {U^{R}}_j & {U^{S}}_j
    \end{bmatrix}}
    \label{eq:low-rank-split}
\end{align}

As shown in Sec.~\ref{sec:construction-and-notation}, the HSS matrix is a multi-level matrix format where each level consists of a single BLR$^2$ matrix. We first introduce the ULV algorithm for the BLR$^2$ matrix in Sec.~\ref{sec:blr2-weak-ulv}. The BLR$^2$-ULV can then be computed at each level of the HSS matrix in order to obtain the HSS-ULV algorithm as shown in Sec.~\ref{sec:ulv-hss-matrix}.

\subsection{Weak admissibility BLR$^2$-ULV factorization}
\label{sec:blr2-weak-ulv}

\begin{figure}
    \centering
    \includegraphics[width=\linewidth]{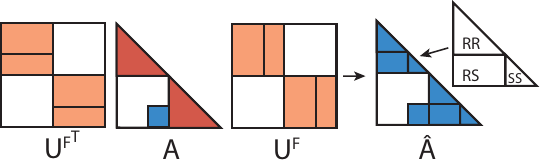}
    \caption{ULV factorization of a weakly admissible BLR$^2$ matrix.}
    \label{fig:blr2-weak-complement-product}
\end{figure}

Alg.~\ref{alg:blr2-weak-ulv} summarizes the BLR$^2$-ULV with weak admissibility. The block diagonal matrix \orangetext{$U^F$} on line 1 is composed of the basis matrices of each row of the BLR$^2$ matrix as shown in Eq.~(\ref{eq:uf-generation}). The resulting product $\hat{A}$ has dense and low rank blocks split into $RR$, $RS$ and $SS$ parts as shown in Eq.~(\ref{eq:dense-split}) and Eq.~(\ref{eq:low-rank-split}), respectively. This is demonstrated in Fig.~\ref{fig:blr2-weak-complement-product} on a 2x2 BLR$^2$ matrix.

\begin{equation}
    \orangetext{U^F} = \orangetext{
        \begin{bmatrix}
            U_0 & 0 \\
            0 & U_1
        \end{bmatrix}
    }
    \label{eq:uf-generation}
\end{equation}

The partial Cholesky factorization on $\hat{A}$ at Line 2 works on the \darkbluetext{$RR$}, \darkbluetext{$RS$} and \darkbluetext{$SS$} blocks of each \darkbluetext{$S_{i,j}$} block. The partial factorization is only performed on the diagonals as a result of the multiplication with the complements in the preceding step
\begin{align}
        L_{i,i}^{RR}{L_{i,i}^{RR}}^{T} \gets Cholesky(\hat{A}_{i,i}^{RR}) \\
        L_{i,i}^{SR} \gets {{L_{i,i}^{RR}}^T}^{-1} \cdot \hat{A}_{i,i}^{SR} \\
        \hat{A}_{i,i}^{SS} \gets \hat{A}_{i,i}^{SS} - L_{i,i}^{SR} \cdot {L_{i,i}^{SR}}^T
    \label{eq:partial-Cholesky}
\end{align}

\begin{figure}
    \centering
    \includegraphics{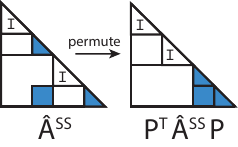}
    \caption{Permutation and Cholesky factorization of partially factorized BLR$^2$ matrix.}
    \label{fig:permute-factorize-blr2-weak}
\end{figure}

Line 3, demonstrated by Fig.~\ref{fig:permute-factorize-blr2-weak}, involves permutation of the partially factorized matrix \darkbluetext{$\hat{A}^{SS}$} which brings all the $SS$ blocks on the lower right corner. This is followed by a dense Cholesky factorization of a smaller matrix of the order of $NB \times rank$. The Cholesky factorization is then performed as follows:
\begin{align}
    \hat{L}^{SS}{\hat{L}^{SS^T}} \gets Cholesky \left(\begin{bmatrix}
        \hat{A}_{00}^{SS} & 0 \\
        \hat{A}_{10}^{SS} & \hat{A}_{11}^{SS}
    \end{bmatrix} \right)
\end{align}

Finally, the matrix \redtext{$A$} is expressed as the following factorization, where \darkbluetext{$\hat{L}$} represents a partially factorized lower diagonal matrix:
\begin{align}
    \redtext{A} = \orangetext{{U^F}^T} \cdot \darkbluetext{\hat{L}} \cdot (P \cdot \darkbluetext{{\hat{L}^{SS}}} \cdot \darkbluetext{{\hat{L}^{SS^T}}} \cdot P^T) \cdot \darkbluetext{\hat{L}^T} \cdot \orangetext{{U^F}}
    \label{eq:blr2-ulv-final-factor}
\end{align}

The solve step for the BLR$^2$-ULV is shown by:
\begin{align}
     x = \orangetext{{U^F}^T} \cdot \darkbluetext{{\hat{L}^{T^{-1}}}} \cdot (P^T \cdot \darkbluetext{{\hat{L}^{SS^{T^{-1}}}}} \cdot \darkbluetext{{\hat{L}^{SS^{-1}}}} \cdot P) \cdot \darkbluetext{\hat{L}^{-1}} \cdot \orangetext{{U^F}} \cdot  b
    \label{eq:blr2-ulv-solve}
\end{align}

\begin{algorithm}[tpb]
  \SetAlgoLined
  \LinesNumbered
  \SetNoFillComment
  \DontPrintSemicolon
  \KwIn{$A$, $U$}
  \tcc{Diagonal product.}
  $\hat{A} \gets {U^F}^T \cdot A \cdot U^F$  \\
  \tcc{Partial factorization.}
  $\hat{A}^{SS} \gets partial\_Cholesky(\hat{A})$  \\
  \tcc{Merge (permute) and factorize.}
  $Cholesky(P^T \cdot \hat{A}^{SS} \cdot P)$
  \caption{ULV factorization of a BLR$^2$ matrix with weak admissibility.}
  \label{alg:blr2-weak-ulv}
\end{algorithm}

The final dense matrix in Alg.~\ref{alg:blr2-weak-ulv} can get large in size for a large rank and problem size. Therefore, even though the leaf level blocks of size $nleaf$ are factorized in $O(N)$, the overall complexity of the algorithm can reach close to $O(N^2)$. The HSS-ULV in the next section shows how the $O(N)$ time complexity can be preserved by exploiting the nested basis.

\subsection{HSS-ULV factorization algorithm}
\label{sec:ulv-hss-matrix}

\begin{algorithm}[tpb]
  \SetAlgoLined
  \LinesNumbered
  \SetNoFillComment
  \DontPrintSemicolon
  \KwIn{$A$, $U$}
    \For{$l\gets max\_level$ \KwTo $1$} { 
      \tcc{Diagonal product.}
      $\hat{A}_{l} \gets {U^F}^T_{l} \cdot A_{l} \cdot U^F_{l}$
      \tcc{Partial factorization.}
      $\hat{A}^{SS}_{l} \gets partial\_Cholesky(\hat{A}_{l})$ \\
      \tcc{Merge (permute).}
      $A_{l-1} \gets P^T_{l} \cdot \hat{A}^{SS}_{l} \cdot P_{l}$
    }
    $Cholesky(A_0)$ \tcc{Final factorization.}
  \caption{ULV factorization of an HSS matrix.}
  \label{alg:hss-ulv}
\end{algorithm}

Alg.~\ref{alg:hss-ulv} describes the ULV factorization of an HSS matrix. The HSS-ULV applies the BLR$^2$-ULV algorithm to each level of the HSS matrix. However, instead of factorizing a dense matrix in the merge step (line 3 in Alg.~\ref{alg:blr2-weak-ulv}), the HSS-ULV algorithm iteratively applies the same procedure to the leftover blocks.

\begin{figure*}[t]
    \centering
    \includegraphics[width=0.7\linewidth]{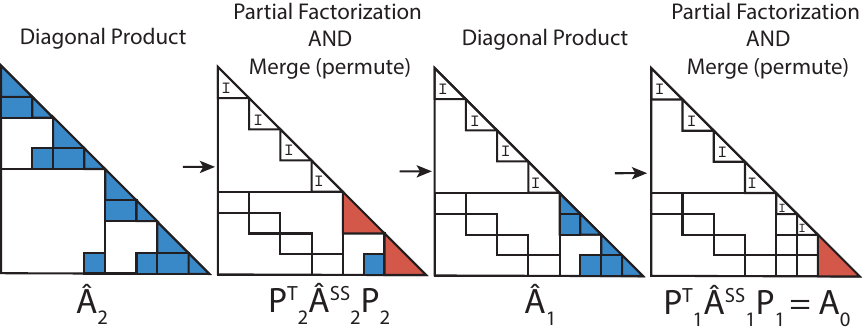}
    \caption{Diagonal product, factorization and permutation steps of the HSS-ULV factorization for a 2 level HSS matrix.}
    \label{fig:hss-ulv}
\end{figure*}

Fig.~\ref{fig:hss-ulv} illustrates the steps taken by the HSS-ULV for 2 iterations of the factorization for a 2 level HSS matrix. A factorization and permutation of the $\hat{A}_2$ matrix leads to the generation of another HSS matrix $P^T_{2} \cdot \hat{A}^{SS}_{2} \cdot P_{2}$, whose diagonal block has a dimension of $2 \times rank$. Another iteration of the ULV factorization of this smaller HSS matrix results in the matrix \redtext{$A_0$} of size $2 \times rank$. Finally a dense Cholesky factorization is performed on \redtext{$A_0$} at line 6 in Alg.~\ref{alg:hss-ulv}. Unlike the merge step of the BLR$^2$-ULV, the final resulting dense matrix for HSS-ULV is much smaller in size. This leads to $O(N)$ time complexity for this algorithm.

The final factorized form of the matrix \redtext{$A$} after the HSS-ULV factorization is very similar to that of the BLR$^2$-ULV in Eq.~(\ref{eq:blr2-ulv-final-factor}). Each lower triangular matrix \darkbluetext{$\hat{L}$} is permuted and multiplied by a \orangetext{$U^F$} corresponding to that level of the matrix. The fully factorized form of the HSS-ULV is shown in Eq.~(\ref{eq:hss-ulv-final-form}).
\begin{align}
    \redtext{A} = &~ \left(\sum_{l=max\_level}^{1}{\orangetext{U^F_l} \cdot \darkbluetext{\hat{L}_l} \cdot P_l}\right) \cdot   \nonumber \\
    & \left( P_0 \cdot \darkbluetext{{\hat{L}^{SS}_0}} \cdot \darkbluetext{{\hat{L}^{SS^T}}_0} \cdot P^T_0\right)  \cdot  \nonumber  \\
    & \left( \sum_{l=1}^{max\_level}P^T_l \cdot \darkbluetext{\hat{L}^T_l} \cdot \orangetext{{U^F_l}^T}\right)
    \label{eq:hss-ulv-final-form}
\end{align}

The solve step of the HSS-ULV is similar to that of BLR$^2$-ULV from Eq.~(\ref{eq:blr2-ulv-solve}). There is again a the introduction of multi-level summation terms similar to the factorization. This solve step is shown below:
\begin{align}
    x = &~\left(\sum_{l=1}^{max\_level}  \orangetext{{{U^F_l}^T}} \cdot \darkbluetext{{\hat{L}_l}^{-1}} \cdot P^T\right) \cdot  \nonumber \\
    & \left( P_0^T \cdot \darkbluetext{{\hat{L}^{S{S^T}^{-1}}_0}} \cdot \darkbluetext{{\hat{L}^{SS^{-1}}}_0} \cdot P_0 \right) \cdot \nonumber \\
    & \left(\sum_{l=max\_level}^{1}{P_l^T  \cdot \darkbluetext{\hat{L}_l^{-1}} \cdot \orangetext{U^F_l}}  \cdot b \right)
\end{align}

\section{Distributed memory execution}
\label{sec:distributed-runtime-systems}

In this section, we elaborate on the distributed memory implementation of \hatrix-DTD, and outline the differences from \strumpack and \lorapo. Sec. \ref{sec:runtime-systems} provides a generic description of runtime systems such as \parsec, followed by an overview of the process distribution strategies followed by \hatrix-DTD, \strumpack and \lorapo in Sec. \ref{sec:process-distribution}. Finally, Sec. \ref{sec:distributed-hss-ulv-parsec} shows how we map the HSS-ULV algorithm to \parsec.

\subsection{Runtime systems.}
\label{sec:runtime-systems}

\begin{figure}
    \centering
    \includegraphics[width=0.8\linewidth]{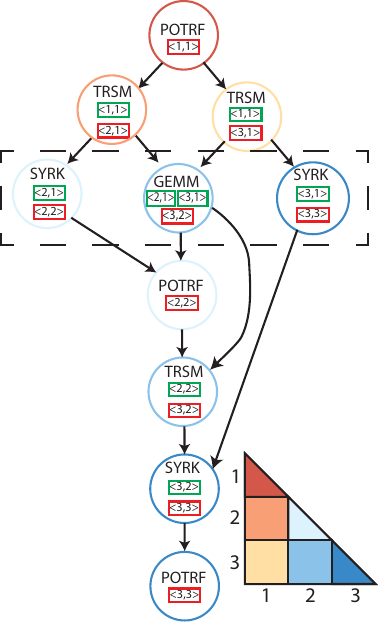}
    \caption{Directed Acyclic Graph (DAG) representation of the block Cholesky factorization of a 3x3 dense matrix. Each node has an associated computation and depends on certain blocks of the matrix (red or green boxes). The red boxes represent RW (Read-Write) dependencies and green boxes represent R (Read) dependencies. The nodes in the dotted blue box shows nodes that can be executed in parallel.}
    \label{fig:dag-Cholesky}
\end{figure}

Fig.~\ref{fig:dag-Cholesky} shows a $3x3$ block Cholesky factorization as a directed acyclic graph (DAG). A runtime system such as \parsec works with such a graph representation of the algorithm in order to factorize the matrix. Each node in the Directed Acyclic Graph in Fig.~\ref{fig:dag-Cholesky} is a `task' with dependencies on some other tasks (nodes) in the graph, shown by the arrows. A task cannot begin execution unless all preceding tasks it depends on have finished their execution. Tasks that do not depend on each other can be executed in parallel. The color of the boundary of each node in the DAG corresponds to the block of the dense matrix that the node updates as a result of the computation taking place in the node.

To illustrate, consider the \texttt{GEMM} task inside the dotted black box. This task has READ dependencies on the $(2,1)$ and $(3,1)$ blocks, and a WRITE dependency on $(3,2)$. The $(2,1)$ and $(3,1)$ blocks that the \texttt{GEMM} must read are WRITE dependencies for the two \texttt{TRSM} tasks that \texttt{GEMM} depends on. Unless both the \texttt{TRSM} tasks before the \texttt{GEMM} finish execution and hand over their blocks to \texttt{GEMM}, it cannot begin execution. Notice that the other two \texttt{SYRK} tasks within the dotted black box also have their respective READ dependencies $(2,1)$ and $(3,1)$ satisfied by the preceding \texttt{TRSM} tasks. Since they have no dependency on the \texttt{GEMM}, the two \texttt{SYRK} tasks and the \texttt{GEMM} inside the dotted black box can be executed in parallel.

\subsection{Distributed HSS-ULV using the \parsec runtime system.}
\label{sec:distributed-hss-ulv-parsec}

\begin{figure*}[htpb]
    \centering
    \includegraphics[width=0.8\linewidth]{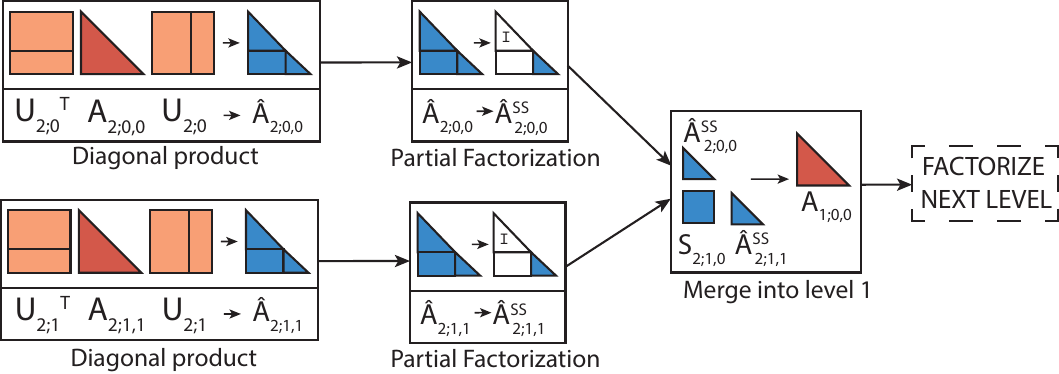}
    \caption{Mapping of the HSS-ULV algorithm to tasks within \parsec for a 2 level HSS matrix. Each block in this diagram represents some computation and its associated dependencies. These tasks represent the factorization of only the first two nodes of the leaf and its subsequent merging into the upper level. Similar steps are followed for other nodes and levels.}
    \label{fig:hss-dag-generation}
\end{figure*}

Fig. \ref{fig:hss-dag-generation} shows the mapping of tasks in \parsec for the HSS-ULV algorithm shown in Alg. \ref{alg:hss-ulv}. We denote the steps shown in Fig.~\ref{fig:hss-ulv} as they are expressed within the tasks of \parsec. The HSS-ULV algorithm operates on the dense, skeleton and bases block matrices of the HSS matrix. These blocks form the dependencies within the tasks. The "Diagonal Product" step on Line 2 of Alg.~\ref{alg:hss-ulv} results in zeroing of the off-diagonal low rank blocks, which makes the partial factorization of each dense block independent of all the other blocks on the same level. This means that the dependencies in the HSS-ULV only come from the merge step on Line 4 of Alg.~\ref{alg:hss-ulv}.

\begin{figure}
    \centering
    \includegraphics[width=0.8\linewidth]{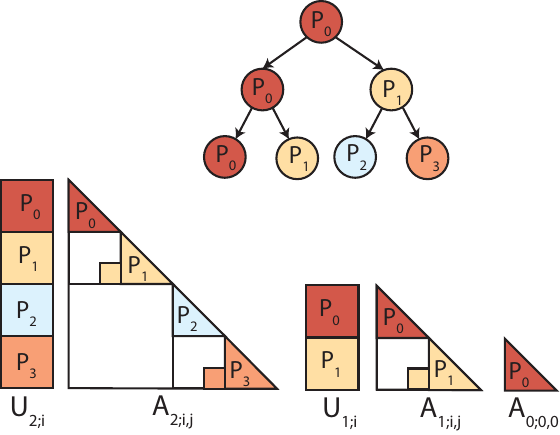}
    \caption{Process distribution used by \hatrix-DTD for a 2 level HSS matrix. Each distinct color represents a separate process. The dense, skeleton and basis blocks are distributed in a row cyclic process distribution at every level.}
    \label{fig:hatrix-hss-process-dist}
\end{figure}

\parsec is able to exploit the inherently parallel factorization of each level. The "Diagonal Product`` and "Partial Factorization`` steps of each dense block on the same level can be executed in an embarrassingly parallel manner. The dependency between the levels exists as a result of the "Merge`` step. As a result of asynchronous execution, the "Merge`` step can begin right after the corresponding partial factorization for its dependencies have finished.

This asynchronous approach is in contrast to \strumpack, where each level executes after the entire previous level has finished factorization. Although the use of SCALAPACK allows \strumpack to perform each level of the HSS-ULV in an embarrassingly parallel manner (assuming different compute resources), the use of fork-join parallelism for the ''Merge`` step means that the parent level cannot begin execution unless the child level has completely finished.

\parsec provides multiple Domain Specific Languages (DSL) to expressing algorithms as DAGs - including Dynamic Task Discovery (DTD) and the Parameterized Task Graph (PTG). The DTD interface is similar to a distributed version of the OpenMP task programming - it submits all tasks on every process, which allows every process to gather all necessary knowledge about the algorithm. However, this means each process discovers the entire task graph, trims the non-local tasks by removing those that do not depend on local tasks, and convert those that depends on communication to and from other processes, and finally execute only the tasks that are local to that process according to their data dependencies. The PTG interface is a custom DSL that allows for an concise parameterized description of the algorithm, and supports an event-driven execution where only local tasks are generated by each process and communications are automatically inferred from the task's dependencies. The PTG interface results in lesser runtime overhead especially when the number of tasks and processes is large as a result of not having to generate the entire task graph on every process. On this paper we focus our efforts into a DTD-based implementation of the algorithm, resulting in \hatrix-DTD.

\subsection{Process distribution.}
\label{sec:process-distribution}

Fig. \ref{fig:hatrix-hss-process-dist} shows the process distribution strategy of \hatrix-DTD for the HSS matrix from Fig. \ref{fig:hss-construction}. Unlike libraries such as SCALAPACK and Elemental~\cite{poulson2013} which make use of block-cyclic and element-cyclic process distribution, respectively, a row-cyclic process distribution is a better fit for HSS-ULV with \parsec. Each block of the HSS matrix that is involved in the HSS-ULV is assigned to a single task. This keeps the number of tasks smaller, which means that the runtime system overhead is better controlled. The blocks owned by $P_0$ and $P_1$ from level 2 are merged into $P_0$ in level 1 as a result of the "Merge`` step in Alg.~\ref{alg:hss-ulv}. Merging blocks into a single process is necessary because of the need to balance the number of tasks with the time consumed by each task. Too many tasks of very little duration will be generated if we proceed with a block-cyclic distribution for the merged block $A_{1;0,0}$.

In contrast to the row cyclic distribution used in \hatrix-DTD, 
\strumpack~\cite{rouetDistributedmemoryPackageDense2015} makes use of a block cyclic distribution for every matrix block of the HSS matrix. This is necessary since \strumpack relies on SCALAPACK for computation on the matrix blocks. Such a process distribution would not be effective in \hatrix-DTD because it would generate too much communication between tasks on the same row (in the block cyclic distribution tasks on the same row will be placed according to the process grid on different processes).

\lorapo~\cite{cao2022a} implements a tile Cholesky algorithm on a block low rank matrix. This means that resolving the off-diagonal triangular solve and trailing sub-matrix updates is the primary bottleneck in the execution of the critical path. A modified block cyclic distribution where each tile is block-cyclically distributed on a process grid is experimentally found to be the best process distribution strategy for \lorapo. Even though \hatrix-DTD makes use of \parsec, such a data distribution in not necessary since the HSS-ULV algorithm ensures that the critical path along the diagonal can be executed in an embarrassingly parallel manner.

\begin{table*}[]
\begin{tabular}{|c|c|c||c|c|c|c|c|c|}
\hline
                                                          & \textbf{\begin{tabular}[c]{@{}c@{}}Construct \\ Max Rank\end{tabular}} & \textbf{\begin{tabular}[c]{@{}c@{}}Leaf \\ Size\end{tabular}} & \textbf{\begin{tabular}[c]{@{}c@{}}Laplace\\ Const. Err.\end{tabular}} & \textbf{\begin{tabular}[c]{@{}c@{}}Laplace\\ Solve Err..\end{tabular}} & \textbf{\begin{tabular}[c]{@{}c@{}}Yukawa\\ Const. Err.\end{tabular}} & \textbf{\begin{tabular}[c]{@{}c@{}}Yukawa\\ Solve Err.\end{tabular}} & \textbf{\begin{tabular}[c]{@{}c@{}}Matern\\ Const. Err.\end{tabular}} & \textbf{\begin{tabular}[c]{@{}c@{}}Matern\\ Solve Err.\end{tabular}} \\ \hline
\multirow{4}{*}{\textbf{HATRIX}}  & 100   & 256   & 1.54e-06  & 4.78e-12  & 2.73e-08  & 3.04e-15  & 9.95e-05    & 3.90e-13    \\ \cline{2-9} 
                                                          & 200  & 256  & 2.89e-07  & 5.48e-12  & 7.63e-09  & 3.50e-15  & 2.34e-05                                                            & 4.51e-12                                                         \\ \cline{2-9} 
                                                          & 200                                                                    & 512                                                           & 1.82e-07                                                           & 6.06e-12                                                           & 5.85e-09                                                          & 3.83e-15                                                         & 1.6e-05                                                               & 4.89e-12                                                         \\ \cline{2-9} 
                                                          & 400                                                                    & 512                                                           & 5.51e-10                                                           & 7.00e-12                                                           & 5.07e-10                                                            & 4.42e-15                                                         & 1.0e-06                                                               & 5.85e-12                                                         \\ \hline
\multirow{4}{*}{\textbf{LORAPO}}                          & 1024                                                                   & 2048                                                          & 1e-8                                                                   & 2.21e-13                                                           & 1e-8                                                                  & 1.33e-13                                                         & 1e-8                                                                  & 1.01e-09                                                         \\ \cline{2-9} 
                                                          & 1500                                                                   & 2048                                                          & 1e-8                                                                   & 2.21e-13                                                           & 1e-8                                                                  & 1.33e-13                                                         & 1e-8                                                                  & 1.01e-09                                                         \\ \cline{2-9} 
                                                          & 1250                                                                   & 4096                                                          & 1e-8                                                                   & 2.21e-13                                                           & 1e-8                                                                  & 1.33e-13                                                         & 1e-8                                                                  & 7.84e-10                                                         \\ \cline{2-9} 
                                                          & 3000                                                                   & 4096                                                          & 1e-8                                                                   & 2.21e-13                                                           & 1e-8                                                                  & 1.33e-13                                                         & 1e-8                                                                  & 7.84e-10                                                         \\ \hline
\multicolumn{1}{|l|}{\multirow{4}{*}{\textbf{STRUMPACK}}} & 100                                                                    & 256                                                           & 1e-8                                                                   & 5.76e-14                                                           & 1e-8                                                                  & 2.13e-15                                                         & 1e-8                                                                  & 1.50e-12                                                         \\ \cline{2-9} 
\multicolumn{1}{|l|}{}                                    & 200                                                                    & 256                                                           & 1e-8                                                                   & 9.05e-11                                                           & 1e-8                                                                  & 1.48e-14                                                         & 1e-8                                                                  & 2.35e-09                                                         \\ \cline{2-9} 
\multicolumn{1}{|l|}{}                                    & 200                                                                    & 512                                                           & 1e-8                                                                   & 3.37e-11                                                           & 1e-8                                                                  & 1.10e-14                                                         & 1e-8                                                                  & 4.44e-10                                                         \\ \cline{2-9} 
\multicolumn{1}{|l|}{}                                    & 400                                                                    & 512                                                           & 1e-8                                                                   & 1.71e-11                                                           & 1e-8                                                                  & 4.04e-14                                                         & 1e-8                                                                  & 9.71e-09                                                         \\ \hline
\end{tabular}
\caption{Impact of rank and kernel for the methods we tested for a constant problem size of 65536.}
\label{tab:rank-accuracy-impact}
\end{table*}

\section{Results}
We run distributed memory tests for 3 implementations of low rank matrix factorization - \hatrix-DTD, \strumpack and \lorapo. Every implementation uses a uniform 2D grid geometry. Distributed memory tests are run on the Fugaku supercomputer at RIKEN, Japan. Each node of Fugaku has a single A64FX CPU with 48 physicals cores divided into 4 NUMA nodes of 12 cores each. Each node has 32 GB of HBM. 

We implemented \hatrix-DTD using the DTD programming interface from \parsec. We make comparisons with \strumpack~\cite{rouetDistributedmemoryPackageDense2015} and \lorapo~\cite{cao2022a} as described in Sec.~\ref{sec:introduction}. We compare three Green's functions from diverse applications such as electrostatics and statistics as shown in Table~\ref{tab:kernels-eval-constants} for each implementation.

\begin{table*}[]
\begin{tabular}{|l|l|l|}
\hline
\textbf{Kernel} & \textbf{Equation} & \textbf{Constants} \\ \hline
 Laplace 2D       &  $f(x,y) = -\ln(\epsilon + dist(x,y))$    &  $\epsilon = 10^{-9}$                  \\ \hline
  Yukawa        &  $f(x,y) = \frac{e^{\alpha \times -(\theta + dist(x,y))}}{(\theta + dist(x,y))}$   & $\alpha = 1$, $\theta = 10^{-9}$                    \\ \hline
  Matern        &  $
    f(x,y) =
    \begin{cases}
    \frac{\sigma^2}{2^{\rho-1} \times \Gamma(\rho)} \times \frac{dist(x,y)}{\mu}^{\sigma} \times K\nu(\sigma, \frac{dist(x,y)}{\mu}),& \text{otherwise} \\
    \sigma^2,              & \text{if } dist(x,y) = 0
    \end{cases}
    $ & $\sigma = 1$, $\mu = 0.03$, $\rho = 0.5$             \\ \hline
\end{tabular}
\caption{Kernels used for evaluation and their constants.}
\label{tab:kernels-eval-constants}
\end{table*}

\subsection{Effect of rank on accuracy}
\label{sec:effect-of-rank-on-acc}

Table~\ref{tab:rank-accuracy-impact} shows the impact of changing the maximum rank and leaf size on the construction and solve error for each tested kernel. The errors are calculated by first generating a normally distributed random vector $b$ and computing the construction error as shown in Eq.~(\ref{eq:construct-error-calc}), and the error of the forward and backward solve as shown in Eq.~(\ref{eq:solver-error-calc}). $A_{dense}$ denotes the full dense matrix, and $A$ denotes the corresponding compressed HSS matrix. The maximum rank of the HSS matrix and the size of the leaf level nodes is capped for all three libraries surveyed. We use the leaf size and maximum rank measurements of solve error from Table~\ref{tab:rank-accuracy-impact} to determine parameters for the weak scaling experiments in Sec. \ref{sec:distributed-mem-weak-scaling} in order to obtain sufficient accuracy for all kernels we benchmark.

\begin{equation}
    err_{construct} = \frac{\norm{A_{dense} \cdot b - A \cdot b}}{\norm{A_{dense} \cdot b}}
    \label{eq:construct-error-calc}
\end{equation}

\begin{equation}
    err_{solve} = \frac{\norm{b - A^{-1} \cdot A \cdot b}}{\norm{b}}
    \label{eq:solver-error-calc}
\end{equation}

The construction error for all cases of \hatrix-DTD decreases as the rank increases, which is to be expected given that a greater rank means a greater number of basis that can be incorporated in the low rank approximation. However, the solve error seems to increase slightly as the rank increases. This slight increase can be attributed to numerical errors. Since \lorapo uses adaptive ranks, specifying the maximum rank leads to choosing ranks that are enough to satisfy the construction error of $10^{-8}$. As a result, increasing the maximum rank from 1024 to 1500 for a leaf size of 2048 does not lead to changing solve error for any kernel. \strumpack allows for a similar tuning of ranks as \hatrix-DTD since both use the HSS matrix. The HSS-ULV algorithm used within \strumpack shows better solve error than \hatrix-DTD for a maximum rank of 100 and leaf size 256 for the laplace 2D kernel. However, the solve error decreases when going from rank 100 to 200 with leaf size 256. The reasons behind this drop in accuracy are unknown.

\begin{figure*}
     \centering
     \begin{subfigure}[b]{0.33\textwidth}
         \centering
         \includegraphics[width=\textwidth]{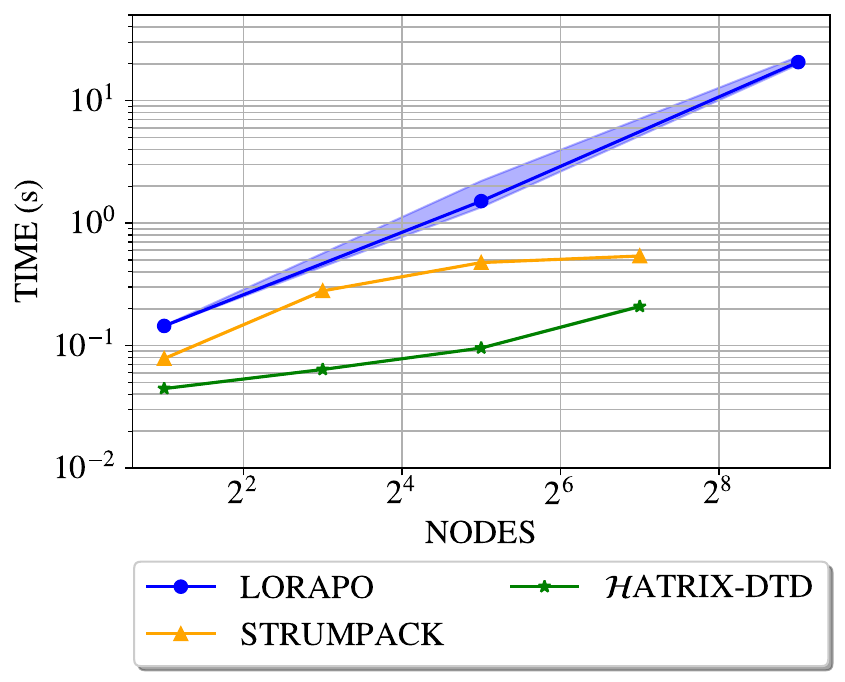}
         \caption{Laplace 2D kernel}
         \label{fig:laplace-2d-weak}
     \end{subfigure}
     \hfill
     \begin{subfigure}[b]{0.33\textwidth}
         \centering
         \includegraphics[width=\textwidth]{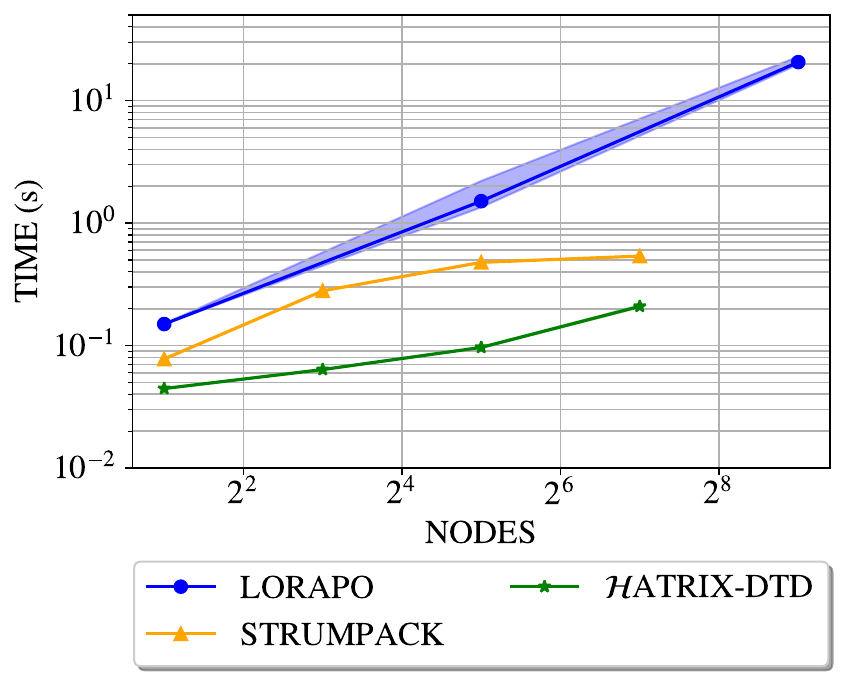}
         \caption{Yukawa kernel}
         \label{fig:yukawa-weak}
     \end{subfigure}
     \hfill
     \begin{subfigure}[b]{0.33\textwidth}
         \centering
         \includegraphics[width=\textwidth]{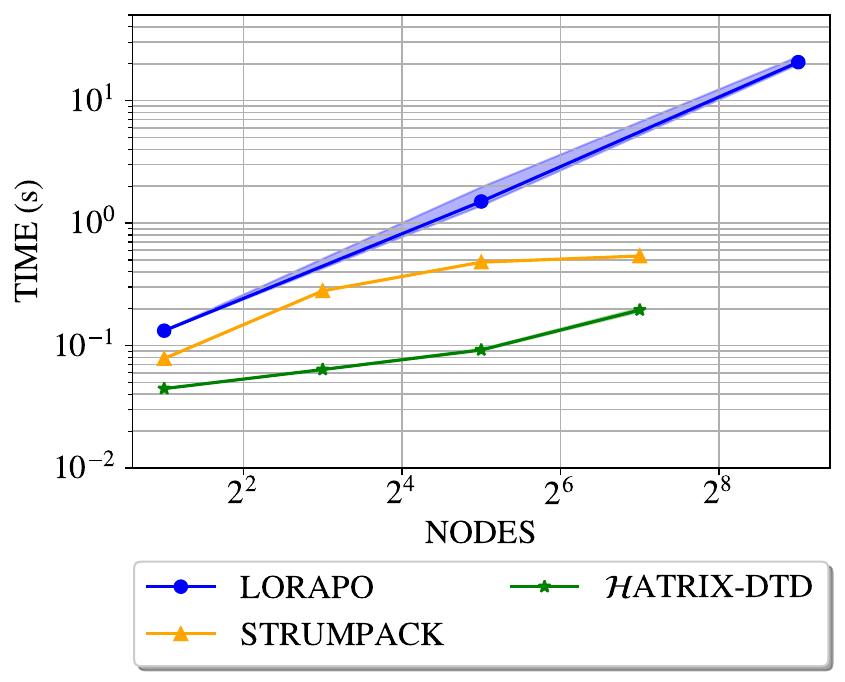}
         \caption{Matern kernel}
         \label{fig:matern-weak}
    \end{subfigure}
    \caption{Weak scaling of factorization time for all the kernels shown in the Table \ref{tab:kernels-eval-constants} for varying problem sizes.}
    \label{fig:factor-time-weak-scaling}
\end{figure*}

\subsection{Distributed memory weak scaling}
\label{sec:distributed-mem-weak-scaling}

Fig.~\ref{fig:factor-time-weak-scaling} shows weak scaling for factorization using \strumpack, \lorapo and \hatrix-DTD on up to 128 nodes of Fugaku. For \hatrix-DTD and \strumpack, the size of the matrix begins at 4096 for 2 nodes and then increases linearly with the number of nodes, until it reaches 262,144 with 128 nodes. The linear increase in problem size and number of processors is done in order to maintain constant work per process, given the $O(N)$ time complexity of the HSS-ULV. The tile Cholesky with the BLR matrix used by \lorapo as shown in Table~\ref{tab:formats-and-algorithms} shows $O(N^2)$ time complexity. Therefore, we start from a problem size of 4096 with 2 nodes and increase the number of nodes by a factor of 16 for every experiment to maintain constant work per node. This means that the problem size reaches 65,536 for 512 nodes. We report the 95\% confidence interval of the mean of the results.

The rank and leaf size are chosen from Sec. \ref{sec:effect-of-rank-on-acc} in order to maintain accuracy that is better than $10^{-11}$ for the laplace 2D kernel, $10^{-14}$ for the Yukawa kernel, and $10^{-9}$ for the matern kernel. We then experiment with combinations of rank and accuracy that provide an acceptable solve error for each problem size and kernel, and show the least time to solution in Fig.~\ref{fig:factor-time-weak-scaling}.

The results in Fig. \ref{fig:factor-time-weak-scaling} show that \hatrix-DTD exhibits better weak scalability than both \strumpack and \lorapo. \lorapo and \hatrix-DTD both make use of the \parsec runtime system, however the tile Cholesky algorithm of \lorapo involves almost $O(N^3)$ communication for the update of the trailing sub-matrix. This, coupled with the fact that the tile Cholesky is constrained by the execution of the critical path of the diagonal add to the poor weak scaling of \lorapo. Further analysis of \lorapo's weak scaling is done in Sec.~\ref{sec:perf-breakdown-lorapo}. \hatrix-DTD and \strumpack both use the HSS-ULV algorithm, however \hatrix-DTD is faster than \strumpack. This is as a result of the asynchronous execution of \parsec, which allows \hatrix-DTD to begin the factorization of the parent level before the entire child level has been factorized. \strumpack, on the other hand, makes use of fork-join parallelism with collective communication, which requires that each level of the HSS matrix be factorized fully before the next level can begin.

\begin{figure*}
     \centering
     \begin{subfigure}[b]{0.33\textwidth}
         \centering
         \includegraphics[width=\textwidth]{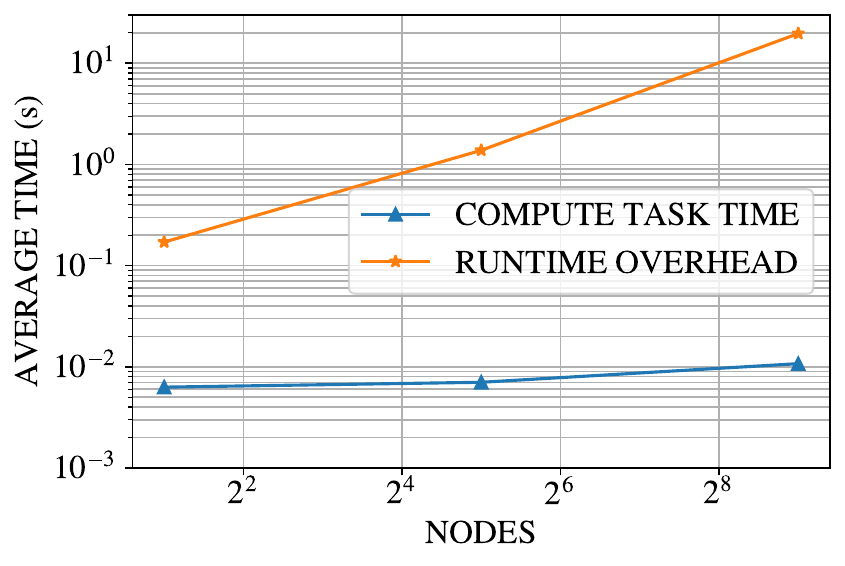}
         \caption{LORAPO}
         \label{fig:lorapo-perf-breakdown}
     \end{subfigure}
     \hfill
     \begin{subfigure}[b]{0.33\textwidth}
         \centering
         \includegraphics[width=\textwidth]{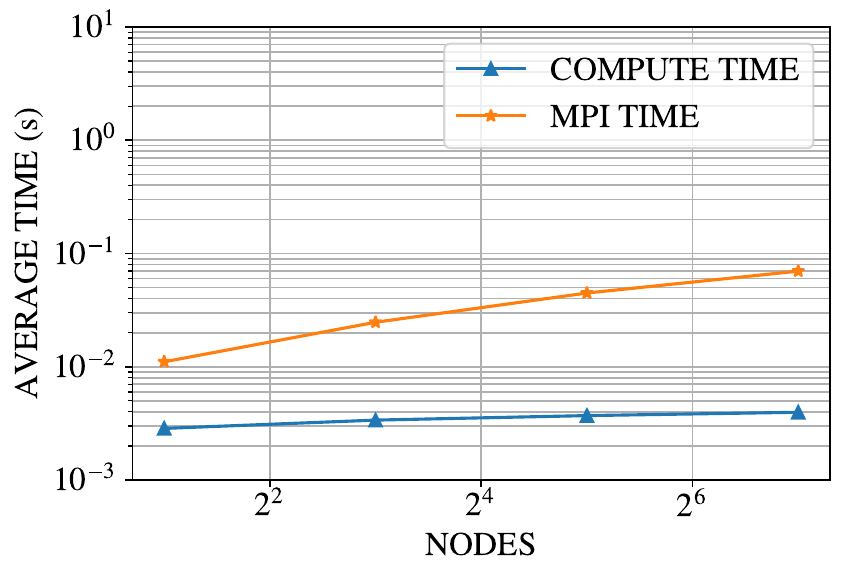}
         \caption{\strumpack}
         \label{fig:strumpack-perf-breakdown}
     \end{subfigure}
     \hfill
     \begin{subfigure}[b]{0.33\textwidth}
         \centering
         \includegraphics[width=\textwidth]{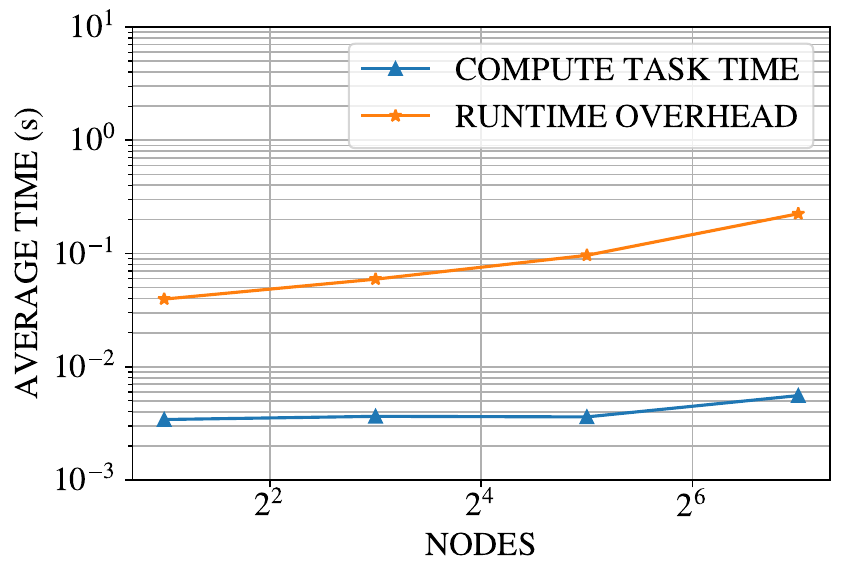}
         \caption{\hatrix-DTD}
         \label{fig:dtd-perf-analysis}
    \end{subfigure}
    \caption{Performance breakdown for the 3 implementations in Fig.~\ref{fig:yukawa-weak}}
    \label{fig:perf-analysis-weak-scaling}
\end{figure*}

\subsection{Performance breakdown of weak scaling}
\label{sec:perf-analysis-weak-scaling}

In this section, we further analyse the reasons behind the weak scaling performance seen in Sec.~\ref{sec:distributed-mem-weak-scaling}. Since all the kernels show similar performance characteristics, we investigate only the Yukawa kernel in further detail.

\subsubsection{Performance breakdown of \lorapo}
\label{sec:perf-breakdown-lorapo}

Fig.~\ref{fig:lorapo-perf-breakdown} shows the performance breakdown for \lorapo~\cite{cao2022a} for the weak scaling graph of the Yukawa kernel shown in Fig.~\ref{fig:yukawa-weak}. We obtain these measurements from the \parsec instrumentation tools that allow for measuring the amount of time that corresponds to time spent inside the actual computational kernels and that for various runtime system management activities.

As pointed out in Sec.~\ref{sec:introduction}, \lorapo uses the tile Cholesky algorithm with the BLR matrix format. The ``COMPUTE TASK TIME'' corresponds to the average time per worker spent inside the actual computational kernels for the Cholesky factorization. The ``RUNTIME OVERHEAD'' corresponds to the average time per worker spent on runtime system management activities such as scheduling, memory management, submitting and executing tasks and deleting previously executed tasks. This also includes various MPI activities such as sending, receiving and polling for messages. The number of workers is the number of physical cores being used for the computation across all the nodes.

It can be seen that overhead of the runtime system far outweighs the amount of time taken for the computation. Moreover, the growth of the overhead is proportional to the time taken for factorization in Fig.~\ref{fig:yukawa-weak}, whereas the growth of ``COMPUTE TASK TIME'' is not. This means that the poor weak scaling of \lorapo can be attributed mainly to the runtime overhead. \lorapo would have had better weak scaling if the runtime overhead would remain constant as the problem size and number of resources is increased.

\subsubsection{Performance breakdown of \strumpack}
\label{sec:perf-breakdown-strumpack}

Fig.~\ref{fig:strumpack-perf-breakdown} shows the breakdown of time that STRUMPACK spends on actual computation vs. MPI for the STRUMPACK weak scaling plot in Fig.~\ref{fig:yukawa-weak}.

The performance statistics in Fig.~\ref{fig:strumpack-perf-breakdown} are obtained from the mpiP tool from LLNL (https://github.com/LLNL/mpiP). The time measurements are averaged over the total number of physical cores used by each experiment. The ``MPI TIME'' shows the time spent by STRUMPACK inside MPI functions such as collective communication. The time does not include the time spent on synchronization. Therefore, the breakdown of computation and communication in Fig.~\ref{fig:strumpack-perf-breakdown} will not add up to the weak scaling performance in Fig.~\ref{fig:yukawa-weak}.  The ``COMPUTE TIME'' shows the time spent on useful computation. The ``COMPUTE TIME'' remains almost the same for every measurement. However, note that the time spent in MPI by each process increases as the number of nodes increases. This means that the MPI communication overhead using the fork-join paradigm leads to inefficient execution in STRUMPACK. Note the "COMPUTE TIME" graph does not show a flat profile in spite of the the HSS-ULV being embarrassingly parallel at each level of the HSS matrix. Apart from the increasing MPI time of the communication, the increasing per process compute time also contributes to the worsening of weak scaling of \strumpack. We show in Sec.~\ref{sec:perf-breakdown-hatrix-dtd} that our implementation in \hatrix-DTD can overcome these limitations.

\subsubsection{Performance breakdown of \hatrix-DTD}
\label{sec:perf-breakdown-hatrix-dtd}

Fig.~\ref{fig:dtd-perf-analysis} shows the performance breakdown of \hatrix-DTD for Fig.~\ref{fig:yukawa-weak}. The measurements are taken in a similar manner to \lorapo in Sec.~\ref{sec:perf-breakdown-lorapo}, i.e. with use of the \parsec instrumentation tools. The ``COMPUTE TASK TIME'' and ``RUNTIME OVERHEAD'' have exactly the same meaning as that of \lorapo in Sec.~\ref{sec:perf-breakdown-lorapo} since both \hatrix-DTD and \lorapo make use of the \parsec runtime system.

Note that the ``COMPUTE TASK TIME'' of \hatrix-DTD is almost completely flat. This means that exactly the same amount of work is being done by each worker when the problem size is increased in proportion to the number of available resources. The final data point shows slightly higher compute time as a result of using a leaf size of 512. This means that the HSS-ULV as implemented in \hatrix-DTD will show perfect weak scaling in the absence of runtime overhead. The increase in run time of \hatrix-DTD in Fig.~\ref{fig:factor-time-weak-scaling} can be attributed to the runtime overhead that shows upward growth as the number of resources is increased. The runtime overhead can be explained by the fact that \parsec's DTD interface generates the entire task graph on every node. This leads to redundant work on each node, which becomes non-trivial as the number of available resources increases.

Note that the compute time per worker for \hatrix-DTD and compute time per thread for \strumpack in Sec. \ref{sec:perf-breakdown-strumpack} are very similar. The runtime overhead of \hatrix-DTD appears higher than that of \strumpack since the MPI barrier and synchronization time is not accounted for in the \strumpack results. However, the \hatrix-DTD overhead shows the entire overhead including synchronization, communication, scheduling and other work by the \parsec runtime system.

\subsection{Increasing problem size with constant resources.}
\label{sec:compare-time-complexity}

\begin{figure}
    \centering
    \includegraphics[width=0.8\linewidth]{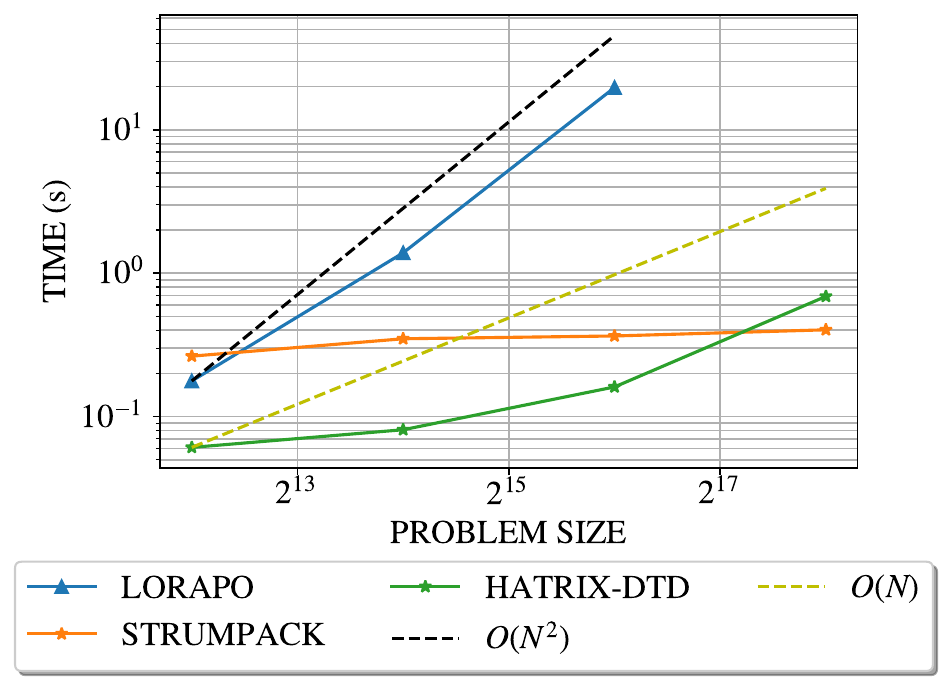}
    \caption{Varying problem sizes with 64 nodes on Fugaku.}
    \label{fig:time-scaling-yukawa}
\end{figure}

Fig.~\ref{fig:time-scaling-yukawa} shows the time taken for factorization for varying problem sizes uses $64$ nodes of Fugaku. \strumpack shows almost uniform time. Since we are using a large number of processes and the computation per process is not very large, the communication time dominates the computation for all cases, and what little computation needs to be done by each process is done in a short time. \lorapo shows $O(N^2)$ scaling up to a problem size of 65,536. \strumpack has an advantage over \hatrix-DTD in this case. Since the runtime overhead in \hatrix-DTD increases as the number of tasks increases, the performance of \hatrix-DTD increases as $O(N)$ even though the amount of computation is small.

\subsection{Impact of leaf size on performance.}
\label{fig:impact-leaf-size}

\begin{figure}
    \centering
    \includegraphics[width=0.8\linewidth]{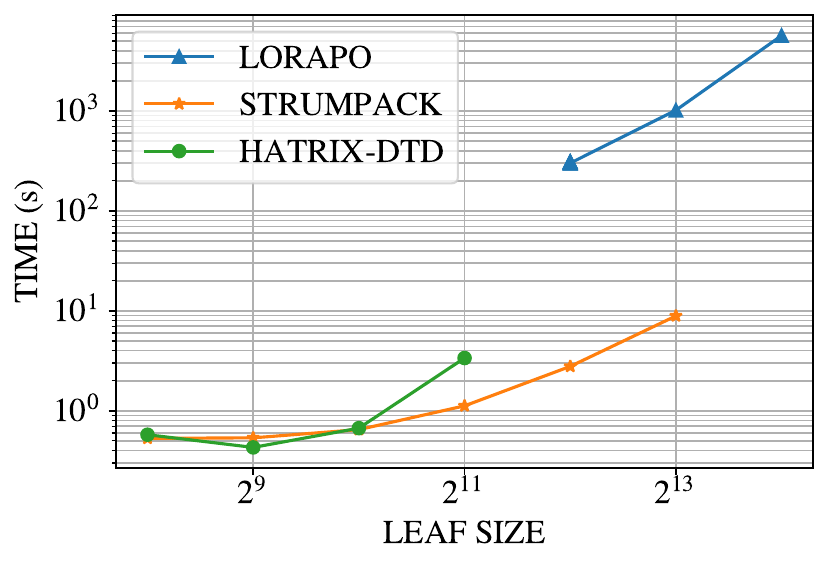}
    \caption{Performance impact of leaf size using 128 nodes and constant problem size of 262,144 for the Yukawa kernel.}
    \label{fig:leaf-size-impact}
\end{figure}

Fig.~\ref{fig:leaf-size-impact} shows the impact on the time for factorization for \hatrix-DTD, \strumpack and \lorapo when using a problem size of 262,144 and 128 nodes of Fugaku. The rank is kept constant at $100$ for \hatrix-DTD and \strumpack and the maximum rank is half the leaf size for \lorapo. The optimal leaf size for \lorapo changes depending on the problem size.  The use of low rank approximation for compressing the dense frontal matrices in the multi-frontal method is an important application of such matrices. The selection of leaf size of the HSS matrix, which correlates to the front size in the multi-frontal solver, is a crucial parameter in justifying the cost of the algorithm. Large leaf sizes can lead to very poor performance of the multi-frontal solver. The fact that \hatrix-DTD is faster than \strumpack when using small leaf sizes shows that \hatrix-DTD can also be used in place of \strumpack to factorize the dense, structured fronts in multi-frontal solvers. Larger leaf sizes for \hatrix-DTD lead to worse performance due to reduction in the amount of available parallelism and more work to do per thread.

\section{Conclusion}

We have proposed a ULV factorization for HSS matrices, and provided an implementation, \hatrix-DTD, using the \parsec runtime system. We have showed that factorization of structured dense matrices arising from a diverse set of Green's functions for a 2D domain can be performed faster using our implementation. This is achieved as a result of the asynchronous runtime system and the lower computational intensity of the HSS-ULV factorization. Using \hatrix-DTD, we show that our implementation has comparable or better accuracy than established state-of-the-art implementations such as \strumpack and \lorapo. Using performance analysis of weak scaling experiments we highlight that our implementation is indeed faster because of a combination of lesser computation and asynchronous resolution of dependencies of the multiple levels of the HSS matrix. As a result of the runtime overhead of \parsec, we have shown that \strumpack can achieve better performance than \hatrix-DTD for large problem size and limited number of nodes.

\begin{acks}
This work was supported by JSPS KAKENHI Grant Number JP20K20624, JP21H03447, JP22H03598. This work is supported by ”Joint Usage/Research Center for Interdisciplinary Large-scale Information Infrastructures” in Japan (Project ID: jh230009-NAHI). This research was also supported by US NSF grant 1909015.
\end{acks}

\bibliographystyle{ACM-Reference-Format}
\bibliography{rioyokotalab}

\end{document}